\documentclass{elsart}
\usepackage{amssymb}

\begin{document}

\begin{frontmatter}

\title{The branching problem in generalized power solutions to differential
equations}

\author{Alejandro S. Jakubi}
\address{Departamento de F\'{\i}sica,
Facultad de Ciencias Exactas y Naturales,
Universidad de Buenos Aires,
Ciudad  Universitaria,  Pabell\'on  I,
1428 Buenos Aires, Argentina.}

\begin{abstract}

Generalized power asymptotic expansions of solutions to differential equations
that depend on parameters are investigated. The changing nature of these
expansions as the parameters of the model cross critical values is discussed.
An algorithm to identify these critical values and generate the generalized
power series for distinct families of solutions is presented, and as an
application the singular behavior of a cosmological model with a
nonlinear dissipative fluid is obtained. This algorithm has been implemented
in the computer algebra system Maple.

\end{abstract}


\end{frontmatter}

\section{Introduction}

Quite frequently physical models depend on parameters and the predictions of
these models, when confronted with observations, yield constraints that fix
the parameters or disqualify the model. Hence, for models described in terms
of a set of differential equations, it becomes of prime importance to
investigate the dependence of solutions  on the parameters. Now, when exact
solutions to these equations are not available, or even when they are
available but they have a complex implicit or parametric form, it may become
useful to obtain asymptotic expansions of these solutions. And when solutions
do not admit power series expansions with integer or rational exponents
(Pusieux series), we are led to consider series with more general terms like
"generalized" powers with real exponents \cite{GS91} or exp-log terms
\cite{vdH}. Though the algebraic issues of generalized power series have already been
investigated (cf. \cite{Loeb}), for applications the main problem
lies in the dependence of the exponents on the parameters, as the ordering of
terms and even the nature of the asymptotic expansion may change as these
parameters reach some critical values.

In this paper, we will discuss the problems involved in the determination of
the critical values of the parameters and an algorithm to calculate the
coefficients and exponents of the generalized power series in the regions of
the parameter space where real solutions admit such an expansion.

Several critical steps in the search for solutions in the form of generalized
power expansions involve calculations with large number of terms and this
number grows very fast with the size of the ODE, making hand calculation
inconvenient. On the other hand, as many of the steps in these calculations
have a systematic nature, use of computer algebra systems appears ideally
suited. For this purpose, we have developed a set of routines in Maple (for a
brief review of a previous implementation, see Ref. \cite{general}).

The plan of this paper is as follows. In section 2, we present the basics of
generalized power expansions of solutions to a class of nonlinear ODES. In
section 3, we review the problem of branching of these solutions as parameters
cross critical values. The algorithm used to find these solutions is presented
in section 4 and it is applied to an example equation arising from a
cosmological model in section 5. Finally, the conclusions are stated in
section 6.

\section{The generalized power expansion}

We will consider
nonlinear ordinary differential equations of the form

\begin{equation} \label{dy}
D[y(t)]=\sum_{i=1}^{N} A_i\, y^{B_i^0}\left(\frac{dy}{dt}\right)^{B_i^1}
\cdots \left(\frac{d^r y}{dt^r}\right)^{B_i^r}=0
\end{equation}

\noindent where the coefficients $A_i$ and the exponent $B_i^0$ may depend on
parameters $p_1,\ldots,p_Q$. In the cosmological setting, $y$ is frequently a
monotonic function of either the scale factor or the Hubble rate in an
expanding universe (cf. \cite{Dav77} \cite{Chi89} \cite{Pav} \cite{AFI}).
Hence $y(t)>0$ and equation (\ref{dy}), if not algebraic, is still well defined.

In case that the general solution to (\ref{dy}) is not available, we are
usually interested in obtaining some information about it in the form of an
asymptotic expansion for the limits $t\to 0^{+}$ and $t\to\infty$, of the form

\begin{equation} \label{yt}
y(t)\sim \sum_{j=1}^\infty c_j t^{n_j}
\end{equation}

\noindent In homogeneous cosmological models, $t$ is the universal time and
these limits frequently correspond to the behavior of the solution near the
initial singularity (the ``Big Bang'') or at large time. Here, $c_j$ and $n_j$
are real, in principle functions of $p_1,\ldots,p_Q$, $c_j\neq 0$ and the
exponents $n_j$ form an ordered set: $n_1< n_2<\cdots$ for $t\to 0^+$ and
$n_1> n_2>\cdots$ for $t\to \infty$. So $c_1 t^{n_1}$ is the leading term,
$n_1\neq 0$, $t^{n_{j+1}}/t^{n_j}\to 0$ in either limit and the set
$\left\{t^{n_j}\right\}$ constitutes an asymptotic scale. Inserting the series
(\ref{yt}) into equation (\ref{dy}), and performing the necessary asymptotic
expansions, we get another generalized power series

\begin{equation} \label{dsum}
D\left[\sum_{j=1}^\infty c_j t^{n_j}\right]\sim \sum_{k=1}^\infty C_k t^{e_k}
\end{equation}

\noindent where $C_k$ and $e_k$ are also real functions of $p_1,\ldots,p_Q$,
and the exponents $e_k$ form an ordered set: $e_1<e_2<\cdots$ for $t\to 0^+$
and $e_1>e_2>\cdots$ for $t\to \infty$. For the equation (\ref{dy}) to admit a
solution with expansion (\ref{yt}), each of the $C_k$ must vanish and this set
of equations fixes in principle the $c_j, n_j$ in (\ref{yt}) up to $r-1$ of
them that remain free and arise from the integration constants of the general
solution of (\ref{dy}) (the arbitrary constant corresponding to the time
translational symmetry is fixed to $0$). Each set of solutions $\{c_j,n_j\}$
yields a series representation of a family of solutions to (\ref{dy}).
Further, the constraints that these $c_j, n_j$ are real and the $n_j$ are
ordered may delimit regions in parameter space where generalized power law
solutions are feasible. If these regions do not contain the values of the
parameters that make physical sense, such solutions have to be discarded even
when they are mathematically correct.

\section{Case branching by parameter variation }

For simplicity, let us consider that equation (\ref{dy}) depends on a single
parameter $p$, so that $n_j=n_j(p)$ and $c_j=c_j(p)$. As the exponents
$\{n_j\}$ have to be ordered, a critical value $p_0$ of the parameter arises
when a pair of consecutive exponents become equal, $n_{j+1}(p_0)=n_j(p_0)$
say. Moving the parameter $p$ across $p_0$ may cause changes in the nature of
the solution, hence in its series expansion, e.g., the asymptotic scale
involved, and these effects show up on the behavior of the function
$\nu_j(p)\equiv n_{j+1}(p)-n_j(p)$ in a neighbourhood of $p_0$. Let us
consider first that $\nu_j(p)$ is analytic at $p_0$ so that, excluding very
special cases, $\nu_j=d_j(p-p_0)+\cdots$ holds with some constant $d_j\neq 0$.
Hence, the terms labeled $j$ and $j+1$ on one side of $p_0$ switch order on
the other side. This effect is of particular importance, as regards to the
analysis of the families of solutions, when the leading behavior changes. For
the series at $p_0$, let us make the expansion of this pair of terms

\begin{equation} \label{ie}
c_jt^{n_j}+c_{j+1}t^{n_{j+1}}=t^{n_j}
\left[c_j+c_{j+1}+c_{j+1}d_j\left(p-p_0\right)\ln t+\cdots\right]
\end{equation}

\noindent If $\lim c_{j+1}(p-p_0)\neq 0$ for $p\to p_0$, as when $c_{j+1}$ is
an arbitrary constant, the series picks a logarithmic term, and this case has
to be dealt with separately. Otherwise, both terms merge into one and a
relabeling of terms occur.

On the other hand, when $p_0$ is a branching point of $\nu_j$ as in
$\nu_j=d_j(p_0-p)^{s/r}+\cdots$, with $s$ odd and $r$ even, an expansion of
the real solution as a generalized power series exists only on one side and
the nature of the solution changes when $p_0$ is crossed.
As an example, let us consider the equation

\begin{equation} \label{1}
\ddot y+y\dot y+\beta y^3 =0
\end{equation}

\noindent
with parameter $\beta$, that arises in several cosmological models (e.g.,
\cite{visco}\cite{Zak}\cite{Reut}\cite{ZPM}), as well as in the analysis of
the Painlev\'e equations \cite{Ince}, and in form invariant equations \cite{Ale}
\cite{Chi97a}. It may be shown that the asymptotic expansion of the general
solution to (\ref{1}) is

\begin{equation} \label{ysing}
y(t)\sim \frac{\alpha}{t} \sum_{n=0}^\infty c_n \gamma^n
 t^{nr}
\end{equation}

\noindent where $\alpha_\pm=[1\pm(1-8\beta)^{1/2}]/(2\beta)$, $r=4-\alpha$ is
the Kowalevski exponent \cite{Yosh}, $r>0$ for $t\to 0^+$ and $r<0$ for
$t\to\infty$,  $c_0=1$, $c_n=c_n(\beta)$ and $\gamma$ is an arbitrary
integration constant. Here, the critical value is $\beta_0=1/8$ and
$\nu_j=r_\pm=\mp 8[2(\beta_0-\beta)]^{1/2}+\cdots$. Hence, real solutions
admit expansion (\ref{ysing}) only for $\beta<1/8$. For $\beta=1/8$, the
general solution has a logarithmic expansion and for $\beta>1/8$, real
solutions become oscillatory.

As the coefficients $c_j$ have to be real and nonvanishing, a critical value
$p_0$ of the parameter arises when it is a root of a coefficient $c_j(p_0)=0$
or  it is a branching point as in $c_j(p)=c_{j0}+d_j(p_0-p)^{s/r}+\cdots$,
with $c_{j0}$ and $d_j$ some constants, $s$ odd and $r$ even. The expansion
(\ref{ysing}) to the solution of equation (\ref{1}) shows the second effect as
its coefficients are proportional to $\alpha$, hence they are complex for
$\beta>1/8$.

\section{The algorithm}

We review briefly the algorithm stated in Ref. \cite{general}.
The objective of the calculations is to obtain a truncation of (\ref{yt}) to
a finite number of terms, say $M$

\begin{equation} \label{yM}
y_M(t)= \sum_{j=1}^M c_j t^{n_j}
\end{equation}

\noindent The method of calculation is iterative, so that constants $c_M,
n_M$, $M>1$, when not free, are determined by
the constants $c_1,n_1,\ldots,c_{M-1},n_{M-1}$ that were found in the
previous steps.

We start with $M=1$, by inserting $y_1=c_1t^{n_1}$ into (\ref{dy}). After
collecting all the terms with the same generalized power, we get a sum of the
form

\begin{equation} \label{dsum1}
D\left[y_1(t)\right]= \sum_{l=1}^R D_l \, t^{f_l}
\end{equation}

\noindent
We note that, in general, the set of the exponents $F=\left\{f_l\right\}_{1\le
l\le R}$ is not ordered. The i-th term of (\ref{dy}) contributes to
(\ref{dsum}) with a term of exponent

\begin{equation} \label{g}
g_i=\sum_{h=0}^r \left(n_1-h\right) B_i^h
\end{equation}

\noindent
and coefficient

\begin{equation} \label{E}
E_i=A_ic_1^{\sum_{h=0}^r B_i^h}\prod_{h=1}^r \left(n_1-h+1\right)^
{\sum_{j=h}^r B_i^j}\equiv \nu_i c_1^{\mu_i}
\end{equation}

As more than one term of (\ref{dy}) may yield the same exponent (balance), we
have $R\le N$. The exponents $g_i$, hence the exponents $f_l$, depend linearly
on $n_1$, and the pair of terms $l$ and $m$ of (\ref{dsum1}), to which
respectively contribute the terms $i(l)$ and $i'(m)$ of (\ref{dy}), balance
for the exponent $n_{1e}^{lm}$

\begin{equation} \label{n1e}
n_{1e}^{lm}=\frac{\sum_{h=0}^r h\left[B_{i(l)}^h-B_{i'(m)}^h\right]}
{\sum_{h=0}^r \left[B_{i(l)}^h-B_{i'(m)}^h\right]}
\end{equation}

We note that both the exponents $g_i$ and the coefficients $E_i$ may also
depend on the parameters $p_1,\ldots,p_Q$ through the exponents $B^0_i$ and
the coefficients $A_i$. In addition, we note that the freedom of multiplying
the equation (\ref{dy}) by integer powers of $y$ and its derivatives (up to
missing solutions where they vanish), that change its exponents by $B^h_i\to
B^h_i+m^h$ for some integers $m^h$, implies that the exponents in (\ref{g}) and
(\ref{E}) are representatives of a class of exponents related by the
transformation laws $g_i\to g_i+\sum_{h=0}^r \left(n_1-h\right) m^h$ and
$\mu_i\to\mu_i+\sum_{h=0}^r m^h$.

If the equation (\ref{dy}) admits a solution of the form (\ref{yt}), and
$y_1(t)$ is its leading term, it holds asymptotically $y(t)\sim y_1(t)$. Hence
on the one hand

\begin{equation}  \label{sim}
D\left[y(t)\right]\sim C_1t^{e_1}\sim D_{l_1}t^{f_{l_1}}
\end{equation}

\noindent
for some index $l_1$ that we identify by sorting the exponents in the set $F$.
For simplicity, only the behavior for $t\to 0^+$ will be considered in the
following. In this case we get $e_1=f_{l_1}\equiv\min(f_l)_{1\le l\le R}$.

On the other hand we require $c_1\neq 0$ and $D_{l_1}=0$. As $f_{l_1}$ is the
minimum of the $\left\{g_i\right\}_{1\le i\le N}$, and we can choose
$\mu_i\neq 0$, at least a pair of terms in (\ref{dy}) must yield this same
exponent. Let us say that this minimum occurs for $i\in I$, so that
$C_1=D_{l_1}=\sum_{i\in I}E_i$ and this coefficient is a function of $c_1$ and
$n_1$.

Then, by solving the set of equations $f_l(n_1)=f_m(n_1)$, $1\le l<m\le R$,
for $n_1$, we obtain the set $N_e$ of equality exponents $n_{1e}$ that make
the exponents $f_l$ exchange order. They delimit intervals within which the
sorting of $F$ has to be carried out separately. Two cases arise: (a) $n_1$
lays inside any of these intervals, (b) $n_1$ is an equality exponent. For
each interval and equality exponent $F$ is sorted and the coefficient of the
term with the minimum exponent $f_{l_1}$ is identified.  In case (b) it is
verified whether $D_{l_1}(c_1)=0$ is satisfied with a nonvanishing root $c_1$,
while in case (a) the equation $D_{l_1}=0$ could determine $n_1$ provided that
a nonvanishing $c_1$ is feasible. Each pair of real numbers $(c_1,n_1)\neq
(0, 0)$ obtained from this analysis corresponds to the leading term of the
expansion of a family of solutions. Thus, subsequent calculations to obtain
higher order terms must proceed separately for each pair. A constant not fixed
by this procedure corresponds in principle to an integration constant of
$y(t)$.

Additional branching occurs due to the dependence on the parameters. As shown
in (\ref{n1e}), the equality exponents depend on $p_1,\ldots,p_Q$ through the
$B_i^0$. If a pair of equality exponents exist, $n_{1e}^{(1)}$ and
$n_{1e}^{(2)}$ say, and they themselves become equal, the equation
$n_{1e}^{(1)}(p_1,\ldots,p_Q)=n_{1e}^{(2)}(p_1,\ldots,p_Q)$ defines a
hypersurface in the parameter space (as the exponents $B_i^0$ may depend on a
subset of the parameters, only a parameter subspace might need to be
considered, and when they depend just on a single parameter, hypersurfaces
become its critical values). Besides, another set of hypersurfaces in the
parameter space might exist where the equality exponents diverge. Hence, the
analysis described before has to be done separately inside each of the
regions of the parameter space delimited by these hypersurfaces, and at their
boundaries.

For $M\ge 2$, when inserted (\ref{yM}) into (\ref{dy}), any factor
$y_M^{(h)B_i^h}$, $h=0,\ldots,r$, with a nonnegative integer exponent $B_i^h$
needs to be expanded, or else expanded asymptotically up to order $M$,
producing a term of the form $K_1t^{(n_1-h)B^h} (1+\cdots+K_Mt^{n_M-n_1})$,
with some constants $K_j$. Then, after such expansions, crossed terms in the
products generate additional terms with larger exponents. For example, when
$M=2$, a term generated by all factors from the leading term except one has an
exponent $g_i'=g_i+n_2-n_1 > g_i$. Thus we see that only the leading term of
(\ref{yt}) can contribute to the leading term of (\ref{dsum}).

Once the solution coefficients $c_1,\ldots, c_{M-1}$ and exponents
$n_1,\ldots,n_{M-1}$ (for a family of solutions) are determined, the main
tasks at step $M$ are:

(i)  Insert $y_M$ in (\ref{dy}) and expand (asymptotically to order $M$) the
powers.

(ii) Collect all the terms with the same power of $t$.

Thus we arrive at an expression of the form

\begin{equation} \label{dsumM}
D\left[y_M(t)\right]= \sum_{l=1}^{R_M} D_l \, t^{f_l}
\end{equation}

\noindent where $D_l$ and $f_l$ are real functions of $p_1,\ldots,p_Q$, and
$R_M\le NM^{r+1}$ as some terms in (\ref{dy}) may balance. In general, the
sequence of exponents $f_1, f_2, \ldots, f_{R_M}$ is not ordered. As $c_M$ and
$n_M$ only appear in $C_k$ for $k\ge M$, the first $M$ terms of (\ref{dsumM}),
once sorted after the order of the exponents, are equal to the first $M$ terms
of expansion (\ref{dsum}). In particular, the first $M-1$ terms are those
already found in step $M-1$ of the iteration. Then, to identify the candidates
for the exponent $e_M$ we

(iii) sort the set of exponents $F_M=\left\{f_l\right\}_{1\le l\le R_M}$, and
pick the $M$-th exponent $f_{l_M}$.

If the equation (\ref{dy}) admits a solution of the form (\ref{yt}), and
$y_M(t)$ is its $M$-term truncation, we get $e_M=f_{l_M}$. The sorting
operation is the most involved part of the whole calculation because of case
branching. As $f_l=f_l(n_M)$, those solution exponents $\{n_{Me}\}$ that make
a pair of exponents equal, $f_l(n_{Me}^{lm})=f_m(n_{Me}^{lm})$ say, and
satisfy $n_{M-1}<n_{Me}$, separate intervals where a given sorting holds. For
$M>1$ these equality exponents $n_{Me}$ arise as solutions to equations of the
form

\begin{equation} \label{lm}
\sum_{j=1}^M \left(\alpha_l^j-\alpha_m^j\right) n_j+
\sum_{i=1}^N\sum_{h=0}^r \left(\beta_l^{ih}-\beta_m^{ih}\right) B_i^h=0
\end{equation}

where the $\alpha_l^j$ are linear functions of the $B_i^h$ (see \cite{general}
for a geometric interpretation of this equation). Furthermore, hypersurfaces
in the parameter space arise as pairs of these equality solution exponents
$n_{Me}$ become equal. In (\ref{lm}) they enter through the $B_i^0$ as well as
through $n_1,\ldots,n_{M-1}$. Hence, step (iii) further divides into:

(iiia) Find the set of the the equality exponents $N_e=\{n_{Me}\}$.

(iiib) Identify the hypersurfaces in the parameter space where the equality
exponents become equal or diverge.

(iiic) Sort the $\{f_l\}$ for each distinct case.

The next steps are:

(iv) Find $n_M$ (if possible) and $C_M$.

(v) Solve $C_M=0$ for either $c_M$ or $n_M$ (if not determined in step (iv)).

A set of routines to deal with steps (iii) and (iv) have been developed in
Maple.

\section{Example}

We will show in this section the application of the algorithm sketched in
section 4 to an equation that depends on several parameters and shows
some of the issues discussed in the previous sections.

In order to treat dissipative processes in cosmology which are not close to
equilibrium, a nonlinear phenomenological generalization of the Israel-Stewart
theory was developed recently \cite{mm}. Scenarios in which this kind of
processes may have occurred include inflation driven by a viscous stress
\cite{mm,m}, and the reheating era at the end of inflation \cite{ZPM,rh}. In a
spatially flat Friedmann-Lemaitre-Robertson-Walker universe, Einstein's
equations together with state and transport equations of the fluid give the
evolution equation for the Hubble rate $H$ \cite{mm}

\begin{eqnarray}\label{14}
&& \left[1-{k^2\over v^2}-\left({2k^2\over3\gamma v^2}\right){\dot{H}
\over H^2}\right]\left\{\ddot{H}+3H\dot{H}+\left(
{1-2\gamma\over\gamma}\right){\dot{H}^2\over H}+{9\over4}
\gamma H^3\right\}  \nonumber\\
&&{}+{3\gamma v^2\over2\alpha}\left[1+\left(
{\alpha k^2\over\gamma v^2}
\right)H^{q-1}\right]H^{2-q}\left(2\dot{H}+3\gamma H^2\right)
-{9\over2}\gamma v^2H^3=0
\end{eqnarray}

\noindent where $\gamma$, $\alpha$, $v$, $q$ and $k$ are parameters describing
the thermodynamical properties of the fluid. We consider an ordinary viscous
fluid so that $1\le\gamma\le 2$, $0<v<1$, $\alpha>0$ and $k>0$. In the
following, we will calculate two term truncations of the generalized power
expansions of the solutions to equation (\ref{14}) in the limit $t\to 0^+$,
corresponding to the behavior of the solutions near the initial singularity.

\subsection{Leading term}

We begin by inserting the leading term $H=c_1t^{n_1}$ into (\ref{14}), looking
for solutions with $n_1<0$ and $c_1>0$. We get the set of exponents

\[
F_1 = \left\{6\,n_1, \,4\,n_1 - 2, \,5\,n_1 - 1, \,3\,
n_1 - 3, \, - n_1\,( q-7), \,6\,n_1 - n_1
\,q - 1\right\}
\]

\noindent
and the set of equality exponents

\[
N_{1e} = \left\{-1, \,{\displaystyle \frac {2}{ q-3}} , \,
{\displaystyle \frac {1}{ q-2}} , \,{\displaystyle \frac {3}{
 q-4}} , \, - {\displaystyle \frac {1}{q}} \right\}
\]

\noindent These, in turn, are function of $q$, and we find that there is a
critical value $q=1$ that make them equal. In effect, this is a distinct value
as all terms of (\ref{dy}) balance and $H=c_1/t$ is an exact solution with
$c_1$ given by the real roots of the cubic equation

$$
\frac{9\gamma}{2}\left [\,{\frac {{\gamma}{v}^{2}}{\alpha}}+\left ({v}^{2}
-\frac{1}{2}\right )\left ({\frac {{k}^{2}}{{v}^{2}}}-1\right )\right ]
{{
 c_1}}^{3}+
$$
$$
3\left [\left (\frac{3}{2}
-v^2\right )\frac{{k}^{2}}{v^2}-1-\,{\frac
{\gamma\,{v}^{2}}{\alpha}}\right ]{{ c_1}}^{2}+
\frac{1}{\gamma}\left (1-3\,{\frac {{k
}^{2}}{{v}^{2}}}\right ){ c_1}+{\frac {2{k}^{2}}{
3{\gamma}^{2}{v}^{2}}}=0
$$

\noindent This shows that $q=1$ delimits different behaviors of the solutions
(cf. \cite{cjmm}). Besides, as an equality exponent diverges for $q=0,2,3,4$,
sorting of $F_1$ must be done separately at these values of $q$ and within the
intervals they delimit, namely $(-\infty,0)$, $(0,1)$, $(1,2)$, $(2,3)$,
$(3,4)$, and $(4,\infty)$. The Table 1 shows the leading exponent for each
case.

\begin{table}
\begin{center}
\begin{tabular}{|c|c|}
\hline
$q<1$ &
\begin{tabular}{cc}
$n_1\le -1$ & $(7-q)n_1$ \\
$-1<n_1\le \displaystyle{\frac{2}{q-3}}$ & $6n_1-n_1q-1$ \\
$n_1>\displaystyle{\frac{2}{q-3}}$ & $3n_1-3$ \\
\end{tabular}
\\ \hline
$1<q\le 3$ &
\begin{tabular}{ccc}
$n_1\le -1$ & ~~ &$6n_1$ \\
$n_1 > -1$ & ~~ &$3n_1-3$ \\
\end{tabular}
\\ \hline
$q>3$ &
\begin{tabular}{cc}
$n_1\le -1$ & $6n_1$ \\
$-1<n_1\le \displaystyle{\frac{2}{q-3}}$ & $3n_1-3$ \\
$n_1>\displaystyle{\frac{2}{q-3}}$ &  $6n_1-n_1q-1$\\
\end{tabular}

\\ \hline
\end{tabular}
\caption[]{Table of leading exponents}
\label{tabla1}
\end{center}
\end{table}

\noindent Thus, we find that the leading exponent switch at $n_1=-1$, as well
as $n_1=2/(q-3)$  for $q<1$ or $q>3$, where terms  balance. Let us start
with these equality exponents. For $n_1=-1$ and $q<1$, we have
$C_1=D(q-7)$ (we denote $D_l$ by $D(f_l)$) where

\[
D(q-7)=\left( - 3\,{\displaystyle \frac {{c_1}^{(6 - q)}\,\gamma }{\alpha
}}  + {\displaystyle \frac {9}{2}} \,{\displaystyle \frac {{
c_1}^{(7 - q)}\,\gamma ^{2}}{\alpha }} \right)\,v^{2}
\]

\noindent so that $c_1=2/(3\gamma)$ is a leading coefficient. For
$q>1$, we have $C_1=D(-6)$ where

$$
D(-6)=
\frac{9\gamma}{2}\left ({v}^{2}
-\frac{1}{2}\right )\left ({\frac {{k}^{2}}{{v}^{2}}}-1\right )
{{
 c_1}}^{3}+
$$
$$
3\left [\left (\frac{3}{2}
-v^2\right )\frac{{k}^{2}}{v^2}-1
\right ]{{ c_1}}^{2}+
\frac{1}{\gamma}\left (1-3\,{\frac {{k
}^{2}}{{v}^{2}}}\right ){ c_1}+{\frac {2{k}^{2}}{
3{\gamma}^{2}{v}^{2}}}
$$

so that we get three leading coefficients

\begin{equation} \label{c1B}
c_1={\frac {2{k}^{2}}{3\gamma\,\left ({k}^{2}-{v}^{2}\right )}}
\,,\quad
{\frac {2}{3\gamma\left (1\pm\sqrt {2}v\right )}}
\end{equation}

provided that $k> v$ in the first case and $v< 1/\sqrt{2}$ in the third
one. For the other equality exponent $n_1=2/(q-3)$ and $q<1$ or $q>3$, we have
$C_1=D(3(5-q)/(q-3))$ where

\[
D\left(3\frac{5-q}{q-3}\right)=
6\,{\displaystyle \frac {{c_1}^{(6 - q)}\,\gamma\,v^{2}}{\alpha \,( q-3)}}  +
\frac{8k^2}{3\gamma v^2(q-3)^3}\left(q-1-\frac{2}{\gamma}\right)c_1^3
\]

\noindent
so that

\begin{equation} \label{c1C}
c_1= \left[{\displaystyle \frac {9\gamma
^{3}\,{ v}^{4}\,(9 - 6\,q + q^{2})}{4{ k
}^{2}\,\alpha \,(  \gamma  - \gamma \,q + 2)}} \right] ^{\frac {1}{ q- 3 }}
\end{equation}

\noindent is the leading coefficient of another family of solutions provided
that $(9 - 6\,q + q^{2})/(  \gamma  - \gamma \,q + 2)>0$. The remaining cases
in Table 1 yield

$$
D((7-q)n_1)=
{\frac {9{{ c_1}}^{7-q}{\gamma}^{2}{v}^{2}}{2\alpha}}
$$

$$
D(6n_1-n_1q-1)=
3\,{\frac {{{ c_1}}^{6-q}\gamma\,{v}^{2}{ n_1}}{\alpha}}
$$

$$
D(3n_1-3)=
\frac{2n_1^2}{3\gamma v^2}
\left(n_1+1-\frac{n_1}{\gamma}\right)
$$

$$
D(6n_1)=
\frac{9\gamma c_1^6}{2}\left(\frac{1}{2}-v^2\right)
\left(1-\frac{k^2}{v^2}\right)
$$

We see that no solution exists in the first two cases, the third case provides
a solution with $n_1=-\gamma/(\gamma-1)$, if $\gamma\neq 1$, and $c_1$
arbitrary, while the fourth case shows that $D(6n_1)=0$ if $k=v$ or
$v=1/\sqrt{2}$. For these values, the calculation of the leading term
has to be done again.

\subsection{Subleading term}

Inserting $H(t)=2/(3\gamma t)+c_2t^{n_2}$ into (\ref{14}) and
expanding the terms with noninteger exponents we get the sets of exponents

$$
F_2 = \{-6,    4{n_2}- 2, q + {n_2}- 6, q + 2{n_2}- 5,  5{n_2}- 1 ,
    3{n_2}- 3, 6{n_2}, {n_2}- 5, 2{n_2}- 4   \}
$$

and the set of equality exponents larger than $-1$, sorted for $q<1$

$$
N_{2e}^<=\left\{ - \frac {1}{2} - \frac {q}{2}, -q\right\}
$$

Thus we find that $e_2=q+n_2-6$ for $n_2\le -q$, while $e_2=-6$ for $n_2>-q$
and $q<1$. In the balancing case $n_2=-q$, we find $C_2=D(-6)$ where

\[
D(-6) =- {\displaystyle \frac {4v^{2}}{3\gamma ^{2
}\,\alpha}}
\left[\alpha  + \gamma (q-2)\left(\frac{3\gamma}{2}\right)^{q}c_2
\right]
\]

so that a solution exists with the subleading coefficient

\begin{equation} \label{c2A}
c_2=\frac {\alpha }{\gamma \,(2 - q)}\left(\frac {2}{3\,\gamma }\right)^{q}
\end{equation}

For the rest of the cases, we have

$$
D(q-6+n_2)=
{\frac {32}{81}}\,{\frac {{3}^{q}{2}^{-q}{\gamma}^{-4+q}\left ({n_2
}+2\right ){v}^{2}{c_2}}{\alpha}}
$$

$$
D(-6)=
-{\frac {32}{81}}\,{\frac {{v}^{2}}{{\gamma}^{5}}}
$$

and neither of them provides a solution.

Following similar steps, inserting $H(t)=c_1/t+c_2t^{n_2}$, with $c_1$ given
by (\ref{c1B}), into (\ref{14}) yields $e_2=n_2-5$ for $n_2\le q-2$, while
$e_2=q-7$ for $n_2>q-2$. At the equality exponent $n_2=q-2$, we obtain the
subleading coefficients corresponding to the three cases of (\ref{c1B})

$$
c_2=-\frac{8k^{6-2q}\left(k^2-v^2\right)^{q-3}v^4}
{9\gamma\alpha q\left(2k^2-v^2\right)}
\left(\frac{3\gamma}{2}\right)^q\,,
$$

$$
\left\{8\sqrt{2}v^4\left(\sqrt{2}v\mp 1\right)^2\left(3\gamma/2\right)^q
\left(1\pm\sqrt{2}v\right)^{q-3}\right\}/
$$
$$
\left\{9\gamma q\alpha\left\{2\sqrt{2}\left[\pm\left(q-1\right)\gamma\mp 2\right]v^4
+2\left[\gamma\left(2k^2-1\right)\left(q-1\right)+4\left(1-k^2\right)
\right]v^3
\right.\right.
$$
$$
+\sqrt{2}\left[\mp\gamma\left(1+2k^2\right)\left(q-1\right)\pm
2\left(4k^2-1\right)\right]v^2
$$
\begin{equation} \label{c2B}
\left.\left.+\left[\gamma\left(1-2k^2\right)\left(q-1\right)-4k^2\right]v\pm
\sqrt{2}\left(q-1\right)k^2\gamma\right\}\right\}
\end{equation}

Inserting $H(t)=c_1t^{2/(q-3)}+c_2t^{n_2}$, with $c_1$ given by (\ref{c1C}),
into (\ref{14}) yields $e_2=(13-3q)/(q-3)+n_2$ for $n_2\le (1+q)/(q-3)$, while
$e_2=2(7-q)/(q-3)$ for $n_2> (1+q)/(q-3)$. At the equality exponent $n_2=
(1+q)/(q-3)$, we obtain the subleading coefficient

\begin{equation} \label{c2C}
c_2=
{\frac {3\left (\left (\left (q-1\right )\gamma-2\right ){v}^{2}+4
\,{k}^{2}\right )\gamma\,\left (q-3\right )}{4{k}^{2}
\left (q-2\right )\left (\left (q-1\right )\gamma-4\right )}}
\left[{\displaystyle \frac {9\gamma
^{3}\,{ v}^{4}\,(9 - 6\,q + q^{2})}{4{ k
}^{2}\,\alpha \,(  \gamma  - \gamma \,q + 2)}} \right] ^{\frac {2}{q - 3}}
\end{equation}

\section{Conclusions}

We have shown some problems that occur when dealing with generalized power
expansions of solutions to nonlinear ordinary differential equations that
depend on parameters, and we have sketched an algorithm that allows
identifying critical values of the parameters and obtaining the series for the
distinct families of solutions by an iterative process.

As an example, we have applied this algorithm to obtain two term truncations of
the series expansions of solutions to a cosmological model filled with a
nonlinear causal viscous fluid. Thus, it is shown as feasible to deal with the
case branching that occurs in series solutions to nonlinear ordinary
differential equations relevant to physical applications.

It deserves to be investigated how the complexity of the algorithm increases
with the order of iteration, and whether this growth puts an effective limit
to practical calculations. In such a case, it would be interesting to know
whether more efficient algorithms could be devised. Also, it would be
interesting to know whether expansions of solutions as shown in this paper
can give information about the integrability of the equation.

\section*{Acknowledgements}
This work was partially supported by the University of Buenos Aires under
Project X223.


\end{document}